\documentclass[10pt]{amsart}
\usepackage{times,amsmath,amsbsy,amssymb,amscd,mathrsfs}
\usepackage{slashbox}\usepackage{graphicx,subfigure,epstopdf,wrapfig,chemarrow}
\usepackage{multicol,multirow}
\usepackage{mathtools}
\usepackage[usenames,dvipsnames,svgnames,table]{xcolor}

\usepackage{comment,enumerate,multicol,xspace}

  \newcounter{mnote}
  \let\oldmarginpar\marginpar
    \renewcommand\marginpar[1]{\-\oldmarginpar[\raggedleft\footnotesize #1]%
    {\raggedright\footnotesize #1}}

\newtheorem{theorem}{Theorem}[section]

\newtheorem{remark}[theorem]{Remark}

\newcommand{\LC}[1]{\textcolor{black}{#1}}
\newcommand{\JH}[1]{\textcolor{black}{#1}}

\begin{document}
\title{Multigrid Methods for A Mixed Finite Element Method of The Darcy-Forchheimer Model} 

\author{Jian Huang Long Chen and Hongxing Rui}

\address[J.~Huang]{School of Mathematics, Shandong University, Jinan, 250100, China.}
\email{yghuangjian@sina.com}

\address[L.~Chen]{Department of Mathematics, University of California at Irvine, Irvine, CA 92697, USA}
\email{chenlong@math.uci.edu}

\address[H.~Rui]{School of Mathematics, Shandong University, Jinan, 250100, China.}
\email{hxrui@sdu.edu.cn}

\thanks{The work of the authors Jian Huang and Hongxing Rui was supported by the National Natural Science Foundation of China Grant No. 11671233, and in part by the Science Challenge Project No. JCKY2016212A502. Long Chen was supported by NSF Grant DMS-1418934, in part by NIH Grant P50GM76516, and in part by the Sea Poly Project of Beijing Overseas Talents. The work of the author Jian Huang was supported by 2014 China Scholarship Council (CSC).}


\begin{abstract}
An efficient \LC{nonlinear} multigrid method for a mixed finite element method of the Darcy-Forchheimer model is constructed in this paper. A Peaceman-Rachford type iteration is used as a smoother to decouple the nonlinearity from the divergence constraint. The nonlinear equation can be solved element-wise with a closed formulae. The linear saddle point system for the constraint is reduced into a symmetric positive definite system of Poisson type. Furthermore an empirical choice of the parameter used in the splitting is proposed and \JH{the resulting multigrid method is robust to the so-called Forchheimer number which controls the strength of the nonlinearity.} By comparing the number of iterations and CPU time of different solvers in several numerical experiments, our multigrid method is shown to convergent with a rate independent of the mesh size and \LC{the Forchheimer number} and with a nearly linear computational cost. 
\end{abstract}

\maketitle

\section{Introduction}
\label{SEC::intro}
Darcy's law
\[\boldsymbol{u} =  - \frac{\boldsymbol{K}}{\mu }\nabla p,\]
\JH{with the permeability tensor $ \boldsymbol{K} $ and the viscosity coefficient $\mu $}, describes the linear relationship between the velocity $\boldsymbol{u}$ of the creep flow and the gradient of the pressure $p$, which is valid when the Darcy velocity $\boldsymbol{u}$ is extremely small~\cite{Aziz}. Forchheimer in~\cite{Forchheimer1901} carried out flow experiments and pointed out that when the velocity is relatively high, Darcy's law should be replaced by the so-called Darcy-Forchheimer (DF) equation by adding a quadratic nonlinear term to the velocity, shown as follows:
\begin{equation}\label{eq:introDF}
\frac{\mu }{\rho }{\boldsymbol{K}^{ - 1}}\boldsymbol{u} + \frac{\beta }{\rho }\left| \boldsymbol{u} \right|\boldsymbol{u} + \nabla p = \boldsymbol{0},
\end{equation}
where $ \rho $  and $ \beta $ represent the density of the fluid and its dynamic viscosity, respectively. \JH{The parameter $ \beta $ is also referred to as the Forchheimer number, which controls the strength of nonlinearity.} A theoretical derivation of the Darcy-Forchheimer equation \eqref{eq:introDF} can be found in~\cite{Ruth1992}. \LC{Equation \eqref{eq:introDF} coupled with the conservation law 
\begin{equation}\label{eq:introdiv}
\text{div}\, \boldsymbol{u} = g
\end{equation}
are usually called Darcy-Forchheimer model.}

In recent years, many numerical methods of the Darcy-Forchheimer model have been developed. 
Girault and Wheeler in~\cite{Wheeler2008} proved the existence and uniqueness of the solution of the
Darcy-Forchheimer model \eqref{eq:introDF}-\eqref{eq:introdiv} by proving the nonlinear operator $\mathcal{A}\left( \boldsymbol{v} \right) = \frac{\mu }{\rho }{\boldsymbol{K}^{ - 1}}\boldsymbol{v} + \frac{\beta }{\rho }\left| \boldsymbol{v} \right|\boldsymbol{v}$ is monotone, coercive and hemi-continuous, and establishing an appropriate inf-sup condition. \LC{Then they considered mixed finite element methods by approximating the velocity and the pressure by piecewise constant and nonconforming Crouzeix-Raviart (CR) elements, respectively. They proved a discrete inf-sup condition and the convergence of the mixed finite element scheme.}
They also proposed a Peaceman-Rachford (PR) type iterative method to solve the discretized nonlinear system and proved convergence of this iterative solver. \JH{In the PR iteration, the nonlinear equation can be decoupled with the divergence constraint and solved in a closed form; see Section \ref{SEC::Nonlinear iteration} for details.} L{\'o}pez, Molina, and Salas in~\cite{Jose2009} carried out numerical tests of the methods proposed in~\cite{Wheeler2008}, and made a comparative study between Newton's method and the PR iterative method. They pointed out that Newton's method is not competitive with the PR iteration. In each iteration, Newton's method needs to evaluate a Jacobian and solves a linear saddle point system, but the PR iteration computes an intermediate solution for a decoupled nonlinear equation and then solves a simplified linear saddle point system. The cost of solving the decoupled nonlinear equation can be negligible in comparison with the Jacobian evaluation. Furthermore the PR iteration required fewer iterations to converge than Newton's method with the same initial guess; see~\cite{Jose2009} for details.

Park in~\cite{Park2005} developed a mixed finite element method with a semi-discrete scheme for the time dependent Darcy-Forchheimer model.
Pan and Rui in~\cite{Pan2012} gave a mixed element method for the Darcy-Forchheimer model based on the Raviart-Thomas (RT) element or the Brezzi-Douglas-Marini (BDM) element approximation \LC{of the velocity and piecewise constant (P0) approximation of the pressure}. 
Rui and Pan in~\cite{Rui2012} proposed a block-centered finite difference method for the Darcy-Forchheimer model, which was thought of as the lowest-order RT-P0 mixed element with proper quadrature formula. Rui, Zhao and Pan in~\cite{RuiZhao2015} presented a block-centered finite difference method for the Darcy-Forchheimer model with variable Forchheimer number $\beta(\boldsymbol{x})$.
Wang and Rui in~\cite{Wang2015} constructed a stabilized CR element for the Darcy-Forchheimer model.
Rui and Liu in~\cite{RuiLiu2015} introduced a two-grid block-centered finite difference method for the Darcy-Forchheimer model. \JH{ Salas, L{\'o}pez, and Molina in~\cite{Jose2013} presented a theoretical study of the mixed finite element method proposed in~\cite{Jose2009}, and showed the well-posedness and convergence.}

\LC{Most of work mentioned above mainly focus on the discretization of the Darcy-Forchheimer model. Except the PR iteration presented in~\cite{Wheeler2008}, no other work concentrates on fast solvers of the discretized nonlinear saddle point system which will be the topic of this paper.}
Multigrid method is one of the most efficient methods on solving the linear and nonlinear elliptic systems. It should be clarified that for nonlinear problems we no longer have a simple linear residual equation, which is the most significant difference between linear and nonlinear systems. The multigrid scheme we used here is the most commonly used nonlinear version of multigrid. It is called the {\em full approximation scheme} (FAS)~\cite{multigrid tutorial} because the problem in the coarse grid is solved for the full approximation rather than the correction; see Section \ref{SEC::Multigrid} for details. 

We shall use piecewise constant (P0) and continuous piecewise linear polynomial (P1) to discretize the velocity and the pressure, respectively. \LC{We refer to~\cite{Jose2013} for the convergence analysis of this scheme and focus on fast solvers in our study.} We shall apply FAS to construct an efficient V-cycle multigrid method for the nonlinear Darcy-Forchheimer model and demonstrate the efficiency of our multigrid method. \LC{Similar application of FAS to a nonlinear saddle point system (for Cahn-Hillard type equations) can be found in~\cite{Wise2010,Wise2013}. \LC{Recall that the success of multigrid method relies on two ingredients: the high frequency can be damped efficiently by the smoother, and the low frequency can be well approximated by the coarse grid correction. Notice that for saddle point systems, both smoothing and coarse grid corrections can easily violate the constraint~\cite{Long2015}. The main difficulty of developing robust and effective multigrid methods for the saddle point system is to design an effective smoother with the consideration of the constraint $\text{div}\, \boldsymbol{u} = g$.}} We shall use the Peaceman-Rachford iteration developed in~\cite{Wheeler2008} as a smoother \LC{since the nonlinearity can be handled efficiently and the constraint is always satisfied after solving a linear saddle point system. To enforce the constraint after the coarse grid correction, we also project the correction into the divergence free subspace. This is in the sprit of the B-S smoother developed in~\cite{Braess1997} for the Stokes equation except here we are dealing with a harder nonlinear equation instead of a linear Stokes equation.}

The most relevant work is~\cite{Jose2009} and our improvement are: 
\begin{enumerate}

\item We reduce the linear saddle point system into a SPD system and demonstrate the efficiency of our approach. 

\item We report a better choice of the splitting parameter $\alpha $ for decoupling the nonlinearity from the constraint rather than the suggested value $ \alpha = 1$ in~\cite{Jose2009} for different values of the Forchheimer number $ \beta $, and show the advantage of our choice by comparing the number of iterations and CPU time. 

\item We carry out some experiments to show the efficiency  of our multigrid solver. \JH{Our method is convergent with a rate independent of the mesh size and the Forchheimer number and with a nearly linear computational cost. Notice that it is not easy to construct a fast solver robust to a critical parameter, see, for example, a linear Stokes-type equation~\cite{Mardal2004,Olshanskii2006}. }

\end{enumerate}

The remainder of this article is organized as follows: The model problem is demonstrated in Section \ref{SEC::The problem}. The mixed weak formulation and the discrete weak formulation are presented in Section \ref{SEC::Weak}. The PR iteration and an efficient solver for the linear saddle point systems are posted in Section \ref{SEC::Nonlinear iteration}. We construct a V-cycle multigrid scheme by applying FAS for the nonlinear problem in Section \ref{SEC::Multigrid}. Some numerical experiments using our multigrid method are carried out in Section \ref{SEC::NES} to verify that the efficiency of our method in comparison with solving this nonliear problem using the other iterative methods. Finally, conclusions, and further ideas are presented in Section \ref{SEC:Conclusions}.

\section{The Problem and Notation}
\label{SEC::The problem}
We consider the steady Darcy-Forchheimer flow of a single phase fluid in a porous medium in a two-dimensional 
bounded domain $ \Omega $, with Lipschitz continuous boundary $\partial \Omega$:
\begin{equation}\label{Eqn::DF equation of the Prob}
\frac{\mu }{\rho }{\boldsymbol{K}^{ - 1}}\boldsymbol{u} + \frac{\beta }{\rho }\left| \boldsymbol{u} \right|\boldsymbol{u} + \nabla p = \boldsymbol{f}\,\,\,\text{in}\,\,\Omega ,
\end{equation}
with the divergence constraint
\begin{equation}\label{Eqn::divergence equation of the Prob}
\text{div}\, \boldsymbol{u} = g\,\,\text{in}\,\,\Omega ,
\end{equation}
and Neumann boundary condition,
\begin{equation}\label{Eqn::N boundary of the Pro}
\boldsymbol{u} \cdot \boldsymbol{n} = {g_N}\,\,\text{on}\,\,{\partial \Omega},
\end{equation}
where $ \boldsymbol{u} $ and $ p $ are the velocity vector and the pressure, respectively; $ \mu $, $ \rho $  and $ \beta $ are given 
positive constants that represent the viscosity of the fluid, its density and its dynamic viscosity, respectively; $ \left| \cdot \right| $ denotes 
the Euclidean vector norm $ \left| \boldsymbol{u} \right| ^2 = \boldsymbol{u} \cdot \boldsymbol{u} $, $ \boldsymbol{n} $ is the unit exterior normal vector to the boundary of the given domain $ \Omega $; $ \boldsymbol{K} $ is the permeability tensor, assumed to be uniformly positive definite and  bounded. According to the divergence theorem, $ g $ and $ g_N $ are given functions satisfying the compatibility condition
\begin{equation}\label{Eqn::compatibility condition}
\int_\Omega  {g\left( \boldsymbol{x} \right)} \, {\rm d}\boldsymbol{x} = \int_{\partial \Omega } {{g_N}\left( \sigma  \right)} \, {\rm d}\sigma .
\end{equation}

We use the standard notation of the Sobolev spaces and the associated norms, see e.g.\cite{sobolev spaces}. 

\section{Weak Formulation}
\label{SEC::Weak}
Following~\cite{Wheeler2008}, we define the function spaces as follows:
\begin{eqnarray*}
X &= & {L^3}{\left( \Omega  \right)^2},\\
M &= & W^{1,{3\over 2}}{\left( \Omega \right)} \cap L_0^2\left( \Omega  \right),
\end{eqnarray*}
where the zero mean value condition
$$L_0^2\left( \Omega  \right) = \left\{ {v \in {L^2}\left( \Omega  \right):\int_\Omega  {v\left( \boldsymbol{x} \right)\, {\rm d}\boldsymbol{x} = 0} } \right\},$$
 is added because $ p $ is only defined by \eqref{Eqn::DF equation of the Prob}-\eqref{Eqn::N boundary of the Pro} up to an additive constant. Given $ \boldsymbol{f} \in {L^3}{\left( \Omega  \right)^2}$, $ g\in {L^{\frac{6}{5}}}{\left( \Omega  \right)}$, and $ g_N \in  {L^{\frac{{3}}{{2}}}}\left( {\partial \Omega } \right) $, the variational formulation of  \eqref{Eqn::DF equation of the Prob}-\eqref{Eqn::N boundary of the Pro} is: find a pair $ \left( \boldsymbol{u},p \right) $ in $ X \times M $ such that
\begin{equation}\label{Eqn::variational formulation1}
\begin{split}
\frac{\mu }{\rho }\int_\Omega  {\left( {{\boldsymbol{K}^{ - 1}}\boldsymbol{u} } \right) \cdot \boldsymbol{\varphi} \, \, {\rm d}\boldsymbol{x}}  &+ \frac{\beta }{\rho }\int_\Omega  {\left| \boldsymbol{u} \right|\left( {\boldsymbol{u} \cdot \boldsymbol{\varphi} } \right)\, {\rm d}\boldsymbol{x}}\\  
&+ \int_\Omega  {\nabla p \cdot \boldsymbol{\varphi} \, \, {\rm d}\boldsymbol{x}}  = \int_\Omega  {\boldsymbol{f} \cdot \boldsymbol{\varphi} \,\, {\rm d}\boldsymbol{x}} ,\quad\forall \boldsymbol{\varphi}  \in X,
\end{split}
\end{equation}
\begin{equation}\label{Eqn::variational formulation2}
\int_\Omega  {\nabla q \cdot \boldsymbol{u}\,\, {\rm d}\boldsymbol{x}}  =  - \int_\Omega  {gq \,\, {\rm d}\boldsymbol{x}}  + \int_{\partial \Omega } {{g_N}q\, \, {\rm d}\boldsymbol{x}} ,\quad\forall q \in M.
\end{equation}

The variational formulation \eqref{Eqn::variational formulation1}-\eqref{Eqn::variational formulation2} and the original problem \eqref{Eqn::DF equation of the Prob}-\eqref{Eqn::N boundary of the Pro} are equivalent by using the Green's formula:
\begin{equation}
\label{Eqn::Green's Formula}
\int_\Omega  {\boldsymbol{v} \cdot \nabla q\,dx}  =  - \int_\Omega  {q \, \mathrm{div} \, \boldsymbol{v}\,dx}  + {\left\langle {q,\boldsymbol{v} \cdot \boldsymbol{n}} \right\rangle _{\partial \Omega }},\quad\forall q \in M,\forall \boldsymbol{v} \in H,
\end{equation}
where 
\[H = \left\{ {\boldsymbol{v} \in {L^3}{{\left( \Omega  \right)}^2}:\mathrm{div}\,\boldsymbol{v} \in {\mathnormal{L}^{\frac{{6}}{{5}}}}\left( \Omega  \right)} \right\}.\]
In~\cite{Wheeler2008}, Girault and Wheeler showed that if the given functions $ g $ and $ g_N $ satisfy the compatibility condition \eqref{Eqn::compatibility condition}, then the problem has a unique solution $ \left( \boldsymbol{u},p \right) $ in $ X \times M $.

Let $\Omega $ be a polygon in two dimensions which can be completely triangulated by triangles. Let $ \mathcal{T}_{\mathnormal{1}} $ be a triangulation of $ \Omega $, and the triangulations $ \mathcal{T}_{\mathnormal{k}} \left( k = 2,3,\ldots\right) $ be obtained form $ \mathcal{T}_{\mathnormal{1}} $ via regular subdivision, i.e. edge midpoints in $ \mathcal{T}_{\mathnormal{k-1}} $ are connected by new edges to form $ \mathcal{T}_{\mathnormal{k}} $. Therefore, $ \mathcal{T}_{\mathnormal{k}} $ is a family of conforming triangulations of $ \overline{\Omega} $,
$$ \overline{\Omega}   = \bigcup\limits_{T \in \mathcal{T}_{\mathnormal{k}}} T\quad \text{for}\;k = 1,2,3, \ldots, $$
The family $ \mathcal{T}_{\mathnormal{k}} $ is shape regular in the sense of Ciarlet~\cite{Ciarlet}. 

We discretize $ \boldsymbol{u} $ and $ p $ in different finite element spaces. The velocity $\boldsymbol{u}$ is approximated in the following space:
\begin{equation}
{X_k} = \left\{ {\boldsymbol{v} \in {L^2}{{\left( \Omega  \right)}^2}:\forall T \in \mathcal{T}_{\mathnormal{k}},\boldsymbol{v}{\vert_{\mathnormal{T}}} \in \mathbb{P}^{\mathrm{2}}_{\mathrm{0}}} \right\},
\label{fes::P0}
\end{equation}
and the pressure $p$ is approximated in the following space:
\begin{equation}
{M_k} = {Q_k} \cap L_0^2\left( \Omega  \right),
\label{fes::P1}
\end{equation}
where $\mathbb{P}_{\mathnormal{m}}$ denotes the space of polynomials of degree $m$, and $ Q_k $ is the linear finite element space
\begin{equation*}
{Q_k} = \left\{ {{q} \in {C^0}{{\left( \bar \Omega  \right)}}:\forall T \in \mathcal{T}_{\mathnormal{k}},\mathnormal{q}{\vert_{\mathnormal{T}}} \in \mathbb{P}_{\mathrm{1}}}\right\}.
\end{equation*}

With these spaces, we can have the $k$-th level discrete formulation of the problem \eqref{Eqn::variational formulation1}-\eqref{Eqn::variational formulation2}:

\begin{equation}\label{Eqn::discrete formulation1}
\begin{split}
\frac{\mu }{\rho }\int_\Omega  {\left( {{\boldsymbol{K}^{ - 1}}\boldsymbol{u}_{k} } \right) \cdot \boldsymbol{\varphi}_{k} \, \, {\rm d}\boldsymbol{x}}  & + \frac{\beta }{\rho }\int_\Omega  {\left| \boldsymbol{u}_{k} \right|\left( {\boldsymbol{u}_{k} \cdot \boldsymbol{\varphi}_{k} } \right)\, {\rm d}\boldsymbol{x}}\\
&+ \sum\limits_{ T \in \mathcal{T}_{\mathnormal{k}}} { \int_T  {\nabla p_{k} \cdot \boldsymbol{\varphi}_{k} \, \, {\rm d}\boldsymbol{x}}}  = \int_\Omega  {\boldsymbol{f} \cdot \boldsymbol{\varphi}_{k} \,\, {\rm d}\boldsymbol{x}} ,\quad\forall \boldsymbol{\varphi}_{k}  \in X_{k},
\end{split}
\end{equation}
\begin{equation}\label{Eqn::discrete formulation2}
\sum\limits_{ T \in \mathcal{T}_{\mathnormal{k}}} {\int_T  {\nabla q_{k} \cdot \boldsymbol{u}_{k}\,\, {\rm d}\boldsymbol{x}}}  =  - \int_\Omega  {gq_{k} \,\, {\rm d}\boldsymbol{x}}  + \int_{\partial \Omega } {{g_N}q_{k}\, \, {\rm d}\boldsymbol{x}} ,\quad\forall q_{k} \in M_{k}.
\end{equation}
By our construction, $$ h_{k - 1} = 2h_{k}, \quad\text{for} \; k = 2,3,\ldots. $$
Note that $ \mathcal{T}_{\mathnormal{k}} $ are nested meshes, and thus
$$  X_{k - 1} \subset X_{k}, \quad M_{k - 1} \subset M_{k}.$$
In~\cite{Jose2013}, the authors demonstrated that the discrete problem \eqref{Eqn::discrete formulation1}-\eqref{Eqn::discrete formulation2} has a unique solution. Moreover, if $\mathcal{T}_{\mathnormal{h}}$ is shape regular with mesh size $h$ and the solution $\boldsymbol{u}$ belongs to ${W^{1,4}}{\left( \Omega  \right)}$ and $ p $ belongs to ${W^{2,\frac{3}{2}}}{\left( \Omega  \right)}$, then the following error estimations are obtained in~\cite[Theorem 4.10]{Jose2013}:
\begin{equation}
\left\| {\boldsymbol{u} - \boldsymbol{u}_{h}} \right\|_{{L^{2}}\left( \Omega  \right)} \le Ch{\left| \boldsymbol{u} \right|_{{W^{1,4}}\left( \Omega  \right)}},
\end{equation}
\begin{equation}
{ {\left\| {\nabla \left( {p - {p_h}} \right)} \right\|} _{{L^{\frac{3}{2}}}\left( T \right)}} \le Ch\left( {{{\left| p \right|}_{{W^{2,\frac{3}{2}}}\left( \Omega  \right)}} + {{\left\| u \right\|}_{{W^{1,4}}\left( \Omega  \right)}}} \right).
\end{equation}

\section{A Nonlinear Iteration}
\label{SEC::Nonlinear iteration}
In this section, we present the Peaceman-Rachford (PR) iterative method developed in~\cite{Wheeler2008} to decouple the nonlinearity and the constraint.

First, choose an initial guess $\left( {\boldsymbol{u}_{k}^{0},{p_k^0}} \right)$ by solving a linear Darcy system:
\begin{equation}
\label{Eqn::darcy1}
\frac{\mu }{\rho }\int_\Omega  {\left( {{\boldsymbol{K}^{ - 1}}\boldsymbol{u}_{k}^0 } \right) \cdot \boldsymbol{\varphi}_{k} \, \, {\rm d}\boldsymbol{x}} + \sum\limits_{ T \in \mathcal{T}_{\mathnormal{k}}} { \int_T  {\nabla p_{k}^0 \cdot \boldsymbol{\varphi}_{k} \, \, {\rm d}\boldsymbol{x}}} = \int_\Omega  {\boldsymbol{f} \cdot \boldsymbol{\varphi}_{k} \,\, {\rm d}\boldsymbol{x}} ,\,\,\forall \boldsymbol{\varphi}_{k}  \in X_{k},
\end{equation}
\begin{equation}
\label{Eqn::darcy2}
\sum\limits_{ T \in \mathcal{T}_{\mathnormal{k}}} {\int_T  {\nabla q_{k} \cdot \boldsymbol{u}_{k}^0\,\, {\rm d}\boldsymbol{x}}}  =  - \int_\Omega  {gq_{k} \,\, {\rm d}\boldsymbol{x}}  + \int_{\partial \Omega } {{g_N}q_{k}\, \, {\rm d}\boldsymbol{x}} ,\,\,\forall q_{k} \in M_{k}.
\end{equation}
The linear Darcy system \eqref{Eqn::darcy1}-\eqref{Eqn::darcy2} can be rewritten in the matrix form as
\begin{equation}
\left[ {\begin{array}{*{20}{c}}
A&B\\
{{B^T}}&0
\end{array}} \right]\left[ {\begin{array}{*{20}{c}}
\boldsymbol{u}\\
p
\end{array}} \right] = \left[ {\begin{array}{*{20}{c}}
\boldsymbol{f_{d}}\\
w
\end{array}} \right],
\end{equation}
where $A$ is the symmetric and positive definite matrix associated to the term 
$$\frac{\mu }{\rho }\int_\Omega  {\left( {{\boldsymbol{K}^{ - 1}}\boldsymbol{u}_{k} } \right) \cdot \boldsymbol{\varphi}_{k} \, \, {\rm d}\boldsymbol{x}},$$
$B$ is the matrix corresponding to
 $$\sum\limits_{ T \in \mathcal{T}_{\mathnormal{k}}} { \int_T  {\nabla p_{k} \cdot \boldsymbol{\varphi}_{k} \, \, {\rm d}\boldsymbol{x}}},$$
and $\boldsymbol{f_{d}}$ and $w$ represent the right hand side of \eqref{Eqn::darcy1} and \eqref{Eqn::darcy2}, respectively.

Then, knowing $\left( {\boldsymbol{u}_{k}^{0},{p_k^0}} \right)$, construct a sequence $\left( {\boldsymbol{u}_{k}^{n + 1},{p_k^{n + 1}}} \right)$ for $n\ge 0$ in two steps.
Let $\alpha$ be a positive parameter chosen to enhance the convergence.

\smallskip
1. A nonlinear step without constraint: knowing $\left( {\boldsymbol{u}_{k}^{n},{p_k^n}} \right)$ compute the intermediate velocity $\boldsymbol{u}_k^{n+\frac{1}{2}}$ by solving the following equation:
\begin{equation}
\label{Eqn::step1 of nonlinear iteration}
\begin{split}
\frac{1}{\alpha }\int_\Omega  {\left( {\boldsymbol{u}_k^{n + \frac{1}{2}} - \boldsymbol{u}_k^n} \right)}  \cdot {\boldsymbol{\varphi} _k}\, {\rm d}\boldsymbol{x} &+ \frac{\beta }{\rho }\int_\Omega  {\left| {\boldsymbol{u}_k^{n + \frac{1}{2}}} \right|} \left( {\boldsymbol{u}_k^{n + \frac{1}{2}} \cdot {\boldsymbol{\varphi} _k}} \right)\, {\rm d}\boldsymbol{x} =  \int_\Omega  {\boldsymbol{f} \cdot \boldsymbol{\varphi}_{k} \,\, {\rm d}\boldsymbol{x}} \\
 - \frac{\mu }{\rho }\int_\Omega  {\left( {{\boldsymbol{K}^{ - 1}}\boldsymbol{u}_k^n} \right)}  \cdot {\boldsymbol{\varphi} _k}\, {\rm d}\boldsymbol{x} &- \sum\limits_{ T \in \mathcal{T}_{\mathnormal{k}}} { \int_T  {\nabla p_{k}^{n} \cdot \boldsymbol{\varphi}_{k} \, \, {\rm d}\boldsymbol{x}}} ,\quad\forall \boldsymbol{\varphi}_{k}  \in X_{k}.
\end{split}
\end{equation}

2. A linear step with constraint: compute $\left( {\boldsymbol{u}_{k}^{n+1},{p_k^{n+1}}} \right)$ with the known $\boldsymbol{u}_k^{n+\frac{1}{2}}$
\begin{equation}
\label{Eqn::step2-1 of nonlinear iteration}
\begin{split}
&\frac{1}{\alpha }\int_\Omega  {\left( {\boldsymbol{u}_k^{n + 1} - \boldsymbol{u}_k^{n + \frac{1}{2}}} \right)}  \cdot {\boldsymbol{\varphi} _k}\, {\rm d}\boldsymbol{x}  + \frac{\mu }{\rho }\int_\Omega  {\left( {{\boldsymbol{K}^{ - 1}}\boldsymbol{u}_k^{n+1}} \right)}  \cdot {\boldsymbol{\varphi} _k}\, {\rm d}\boldsymbol{x} + \sum\limits_{ T \in \mathcal{T}_{\mathnormal{k}}} { \int_T  {\nabla p_{k}^{n+1} \cdot \boldsymbol{\varphi}_{k} \, \, {\rm d}\boldsymbol{x}}} \\
 &= \int_\Omega  {\boldsymbol{f} \cdot \boldsymbol{\varphi}_{k} \,\, {\rm d}\boldsymbol{x}} 
 - \frac{\beta }{\rho }\int_\Omega  {\left| {\boldsymbol{u}_k^{n + \frac{1}{2}}} \right|} \left( {\boldsymbol{u}_k^{n + \frac{1}{2}} \cdot {\boldsymbol{\varphi} _k}} \right)\, {\rm d}\boldsymbol{x}  ,\quad\forall \boldsymbol{\varphi}_{k}  \in X_{k},
\end{split}
\end{equation}
\begin{equation}
\label{Eqn::step2-2 of nonlinear iteration}
\sum\limits_{ T \in \mathcal{T}_{\mathnormal{k}}} {\int_T  {\nabla q_{k} \cdot \boldsymbol{u}_{k}^{n+1}\,\, {\rm d}\boldsymbol{x}}}  =  - \int_\Omega  {gq_{k} \,\, {\rm d}\boldsymbol{x}}  + \int_{\partial \Omega } {{g_N}q_{k}\, \, {\rm d}\boldsymbol{x}} ,\quad\forall q_{k} \in M_{k}.
\end{equation}

A key observation in~\cite{Wheeler2008} is that because the test functions ${\boldsymbol{\varphi} _k}$, the solution $ \boldsymbol{u}_k^{n + \frac{1}{2}} $, and $ \nabla p_{k}^{n} $ are constant in each element $ T $, the nonlinear step \eqref{Eqn::step1 of nonlinear iteration} can be solved in a closed-form:
\begin{equation}
\label{Eqn::explicit solution of nonlinear term}
 \boldsymbol{u}_T^{n + \frac{1}{2}} = \frac{1}{\gamma}  \boldsymbol{F}_T^{n + \frac{1}{2}}
\end{equation} 
where 
\begin{align*}
 \boldsymbol{F}_T^{n + \frac{1}{2}} & = \frac{1}{\alpha}\boldsymbol{u}_T^{n} - \frac{\mu}{\rho}\boldsymbol{K}_T^{-1}\boldsymbol{u}_T^{n} - \nabla_T p_{k}^{n} + \boldsymbol{f}_T,  \\
  \boldsymbol{K}_T^{-1} & = \frac{1}{\left| {T} \right| } \int_T {\boldsymbol{K}^{-1}\left(\boldsymbol{x}\right) \, {\rm d}\boldsymbol{x}}, \\
  \gamma & = \frac{1}{2\alpha} + \frac{1}{2}\sqrt{\frac{1}{\alpha^2} + 4\frac{\beta}{\rho}\left| {\boldsymbol{F}_T^{n + \frac{1}{2}}} \right|}.
\end{align*}

\par

In the second step, the linear system \eqref{Eqn::step2-1 of nonlinear iteration}-\eqref{Eqn::step2-2 of nonlinear iteration} can be rewritten in the following matrix form:
\begin{equation}
\label{Eqn::matrix form of nonlinear iteration}
\left[ {\begin{array}{*{20}{c}}
{A_{\alpha} }&B\\
{{B^T}}&0
\end{array}} \right]\left[ {\begin{array}{*{20}{c}}
\boldsymbol{u}\\
p
\end{array}} \right] = \left[ {\begin{array}{*{20}{c}}
\boldsymbol{f_{n + \frac{1}{2}}}\\
w
\end{array}} \right],
\end{equation}
where $A_{\alpha}$ is the matrix corresponding to the bilinear form
$$ \frac{1}{\alpha }\int_\Omega  {\left( {\boldsymbol{u}_k^{n + 1}} \right)}  \cdot {\boldsymbol{\varphi} _k}\, {\rm d}\boldsymbol{x} + \frac{\mu }{\rho }\int_\Omega  {\left( {{\boldsymbol{K}^{ - 1}}\boldsymbol{u}_k^{n + 1}} \right)}  \cdot {\boldsymbol{\varphi} _k}\, {\rm d}\boldsymbol{x}, $$
and $\boldsymbol{f_{n + \frac{1}{2}}}$ is the vector corresponding to 
$$ \int_\Omega  {\boldsymbol{f} \cdot \boldsymbol{\varphi}_{k} \,\, {\rm d}\boldsymbol{x}} + \frac{1}{\alpha }\int_\Omega  {\left( {\boldsymbol{u}_k^{n + \frac{1}{2}}} \right)}  \cdot {\boldsymbol{\varphi} _k}\, {\rm d}\boldsymbol{x} - \frac{\beta }{\rho }\int_\Omega  {\left| {\boldsymbol{u}_k^{n + \frac{1}{2}}} \right|} \left( {\boldsymbol{u}_k^{n + \frac{1}{2}} \cdot {\boldsymbol{\varphi} _k}} \right)\, {\rm d}\boldsymbol{x}.$$ 

In~\cite{Wheeler2008}, the authors proved that \eqref{Eqn::darcy1}-\eqref{Eqn::darcy2} and \eqref{Eqn::step2-1 of nonlinear iteration}-\eqref{Eqn::step2-2 of nonlinear iteration} have a unique solution. The PR iterative method is convergent for an arbitrary choice of the initial guess $\left( {\boldsymbol{u}_{k}^{0},{p_k^0}} \right)$ and an arbitrary positive $ \alpha $. Numerically, different choices of $\alpha$ will affect the convergence rate of the nonlinear iteration. We shall report a choice of $\alpha$ in Section \ref{SEC::NES}. 

We can reduce the linear saddle point system into a SPD system when we implement the PR iteration. Because of $ A $ and $ A_{\alpha} $ are symmetric positive definite operators, without loss of generality, we take \eqref{Eqn::matrix form of nonlinear iteration} as an example to expound an idea as follows.

Eliminate $ \boldsymbol{u} $ from the first equation of \eqref{Eqn::matrix form of nonlinear iteration}, i.e.
\begin{equation}
\label{Eqn::Solve u of the SPD system}
\boldsymbol{u} = A_{\alpha}^{-1}\left(\boldsymbol{f_{n + \frac{1}{2}}} - Bp \right),
\end{equation}
and then, substituting to the second equation of \eqref{Eqn::matrix form of nonlinear iteration}, we get
\begin{equation}
Mp = b,\label{Eqn::SPD system}
\end{equation}
where $ M = B^{T}A_{\alpha}^{-1}B, b = B^{T}A_{\alpha}^{-1}\boldsymbol{f_{n + \frac{1}{2}}} - w $. After solving \eqref{Eqn::SPD system}, we can get $ \boldsymbol{u} $ by solving \eqref{Eqn::Solve u of the SPD system}.

Since $ A_{\alpha} $ is block-diagonal, $ A_{\alpha}^{-1} $ can be formed easily. Indeed equation \eqref{Eqn::SPD system} is the linear finite element discretization of an elliptic equation in the primary formulation. The equivalence between \eqref{Eqn::Solve u of the SPD system}-\eqref{Eqn::SPD system} and \eqref{Eqn::matrix form of nonlinear iteration} is obvious. Solving the SPD system \eqref{Eqn::SPD system} is much easier than the saddle point system \eqref{Eqn::matrix form of nonlinear iteration} and many fast solvers are available. In our numerical experiments, we use the direct solver built in ${\rm MATLAB}^{\copyright}$ to solve \eqref{Eqn::SPD system}. We could also use the multigrid solver, but due to the relative-small size of the linear SPD system we have tested, the direct solver is faster.

In the continuous level, the Darcy-Forchheimer equation can be rewritten into a nonlinear primary formulation. For simplicity, we assume that the permeability is a scalar. Taking the norm of equation \eqref{Eqn::DF equation of the Prob}, we obtain
\begin{eqnarray*}
\frac{\beta }{\rho }{\left| \boldsymbol{u} \right|}^2 + \frac{\mu }{\rho K}\left| \boldsymbol{u} \right| - \left| \nabla p - \boldsymbol{f}\right| = 0,
\end{eqnarray*}
and can solve for $\left| \boldsymbol{u} \right|$
\begin{equation*}
\left| \boldsymbol{u} \right| = \frac{-\frac{\mu }{\rho K} + \sqrt{\left({\frac{\mu }{\rho K}}\right)^2 + 4\frac{\beta }{\rho }\left| \nabla p - \boldsymbol{f}\right|}}{2\frac{\beta }{\rho }}.
\end{equation*} 
and consequently $\boldsymbol{u}$ 
\begin{equation*}
\boldsymbol{u} = -\frac{\nabla p - \boldsymbol{f}}{\frac{\mu }{\rho K} + \frac{\beta }{\rho }\left| \boldsymbol{u} \right|} = -\frac{2\left(\nabla p - \boldsymbol{f}\right)}{\frac{\mu }{\rho K} + \sqrt{\left({\frac{\mu }{\rho K}}\right)^2 + 4\frac{\beta }{\rho }\left| \nabla p - \boldsymbol{f}\right|}}.
\end{equation*}
Then substituting back to \eqref{Eqn::divergence equation of the Prob}, we get the primary formulation of pressure $p$ only
\begin{equation}\label{eq:primary}
-\nabla \cdot \left(\frac{2\left(\nabla p - \boldsymbol{f}\right)}{\frac{\mu }{\rho K} + \sqrt{\left({\frac{\mu }{\rho K}}\right)^2 + 4\frac{\beta }{\rho }\left| \nabla p - \boldsymbol{f}\right|}}\right) = g.
\end{equation}
Its well-posedness can be found in~\cite{Pan2012}.

In the discretization level, we could also eliminate the piecewise constant velocity and obtain an equivalent $P_1$ discretization of \eqref{eq:primary}. 
However, we only eliminate $ \boldsymbol{u} $ of the linear system \eqref{Eqn::matrix form of nonlinear iteration} in the PR iteration rather than that of the nonlinear equation \eqref{Eqn::discrete formulation1} because we still need to solve the resulting nonlinear equation. The PR iteration corresponds to a variant of Picard iteration for solving \eqref{eq:primary}.
We stick to the mixed formulation as the convergence of the PR iteration has been rigorously proved in~\cite{Wheeler2008}.

\section{Non-linear Multigrid Algorithm}
\label{SEC::Multigrid}
In this section, we consider a generic system of nonlinear equations,
$$ \mathcal{L}\left( \boldsymbol{z} \right) = \boldsymbol{s} $$
where $ \boldsymbol{z}, \boldsymbol{s} \in \boldsymbol{R^{n}}$. Suppose that $ \boldsymbol{v} $ is an approximation to the exact solution $ \boldsymbol{z} $. Define the error $\boldsymbol{e} $ and the residual $ \boldsymbol{r} $:
\begin{align*}
\boldsymbol{e} & = \boldsymbol{z} - \boldsymbol{v}, \\
\boldsymbol{r} & = \boldsymbol{s} - \mathcal{L}\left(\boldsymbol{v}\right).
\end{align*}
Quantities in the $k$-th level will be denoted by a subscript $k$. 
  
\par
Because of the iterative nature, multigrid ideas should be effective on the nonlinear problem. The multigrid scheme here we used for this nonlinear problem is the most commonly used nonlinear version of multigrid. It is called the {\em full approximation scheme} (FAS)~\cite{multigrid tutorial} because the problem in the coarse grid is solved for the full approximation $ {\boldsymbol{z}_{k-1}} = I_{k}^{k-1}\boldsymbol{v}_{k} + \boldsymbol{e}_{k-1} $ rather than the error $\boldsymbol{e}_{k-1}$. A two-level FAS is described as follows. 

~\\
\textbf{Full Approximation Scheme (FAS).}
~\\
\begin{enumerate}

\item Pre-smoothing:  For $ 1 \leq j \leq m, $ relax $ m $ times with an initial guess $ \boldsymbol{v}^{0} $ by  $ {\boldsymbol{v}^j} = {R_k}{\boldsymbol{v}^{j - 1}} $. The current approximation $ {\boldsymbol{v}_k} = {\boldsymbol{v}^m} $.
~\\
\item Restrict the current approximation and its fine grid residual to the coarse grid: $ \boldsymbol{r}_{k-1} = I_{k}^{k-1}\left(\boldsymbol{s}_{k} - \mathcal{L}_{k}\left( \boldsymbol{v}_{k}\right)\right) $ and $ \boldsymbol{v}_{k-1} = I_{k}^{k-1}\boldsymbol{v}_{k} $.
~\\
\item Solve the coarse grid problem: $ \mathcal{L}_{k-1}\left(\boldsymbol{z}_{k-1}\right) = \mathcal{L}_{k-1}\left(\boldsymbol{v}_{k-1}\right) + \boldsymbol{r}_{k-1} $.
~\\
\item Compute the coarse grid approximation to the error: $ \boldsymbol{e}_{k-1} = \boldsymbol{z}_{k-1} - \boldsymbol{v}_{k-1} $.
~\\
\item Interpolate the error approximation up to the fine grid and correct the current fine grid approximation: $ \boldsymbol{v}_{m + 1} \leftarrow \boldsymbol{v}_{k} + I_{k-1}^{k}\boldsymbol{e}_{k-1}$.
~\\
\item Post-smoothing: For $ {m + 2} \leq j \leq {2m + 1}, $ relax $ m $ times by  $ {\boldsymbol{v}^j} = {R_k}'{\boldsymbol{v}^{j - 1}} $.
\end{enumerate}
then we get the approximate solution $ {\boldsymbol{v}^{2m + 1}} $.  Here $ m $ denotes the number of pre-smoothing and post-smoothing steps, $ {R_k} $ denotes the chosen relaxation method, and $ I_{k}^{k-1} $ is an intergrid transfer operator from the fine grid to the coarse grid. \LC{As usual, the V-cycle will be obtained by applying the two-level FAS to the solve the nonlinear equation in Step 3.}
\par

We choose the PR iteration \eqref{Eqn::step1 of nonlinear iteration}-\eqref{Eqn::step2-2 of nonlinear iteration} as the smoother ${R_k}$ and the nonlinear solver in the coarsest grid. We switch the ordering of the linear and nonlinear steps \JH{of the PR iteration} in the post-smoothing step in order to keep the symmetry of the V-cycle. It is worth pointing out that although the chosen finite element spaces are nested, the constrained subspaces are non-nested when we interpolated the correction of the velocity, which was obtained in the coarser space, to the finer space. Namely, if we directly interpolated the correction obtained on the coarser grid to the finer grid, the approximation we got may not satisfy the divergence equation in this Darcy-Forchheimer model. Therefore we construct a $ L^2 $ projection to map the correction obtained before into the constrained space in the fine grid which can be realized by solving a saddle point system:
\begin{equation}
\label{Eqn::L^2 projection}
\left[ {\begin{array}{*{20}{c}}
{A_{\boldsymbol{\delta}} }&B\\
{{B^T}}&0
\end{array}} \right]\left[ {\begin{array}{*{20}{c}}
\boldsymbol{\delta}\\
\theta
\end{array}} \right] = \left[ {\begin{array}{*{20}{c}}
\boldsymbol{0}\\
B^T\boldsymbol{e_{u}}
\end{array}} \right],
\end{equation}
where $ A_{\boldsymbol{\delta}} $ is the matrix corresponding to 
$$\frac{\mu }{\rho }\int_\Omega  {\left( {{\boldsymbol{K}^{ - 1}}\boldsymbol{\delta} } \right) \cdot \boldsymbol{\varphi}_{k} \, \, {\rm d}\boldsymbol{x}} + \frac{\beta }{\rho }\int_\Omega  {\left| \boldsymbol{\delta} \right|\left( {\boldsymbol{\delta} \cdot \boldsymbol{\varphi}_{k} } \right)\, {\rm d}\boldsymbol{x}},$$
$ \boldsymbol{\delta}, \theta $ represent the error between the restriction of the approximation of velocity and pressure on the finer grid and their approximation obtained on the coarser grid, respectively, and $ \boldsymbol{e_{u}} $ is the prolonged correction to the fine space. For non-nested constrained subspaces, an additional projector is usually needed to preserve the constraint~\cite{Braess1997}.

Again, \eqref{Eqn::L^2 projection} can be reduced to a SPD system. We can get $ \boldsymbol{\delta} = A_{\boldsymbol{\delta}}^{-1}B\theta $ through the idea demonstrated in Section \ref{SEC::Nonlinear iteration}. Then we obtain a corrected approximation of velocity $ \boldsymbol{v} = \boldsymbol{v} - \boldsymbol{\delta}$, which satisfies the divergence equation.

\begin{remark}
\LC{When RT or BDM element is used to discretize the velocity and the pressure is piecewise constant, we may use patch-wise smoothers designed for $H({\rm div})$ problems; see~\cite{Arnold1997,Arnold2000}. The constraint can be preserved in these smoothers. A rigorous proof for the convergence of a multigrid method using constrained smoothers for linear saddle point systems can be found in~\cite{Long2015,LongMathComp}. Note that in this paper, we consider continuous pressure discretization and nonlinear saddle point systems and thus neither the constrained smoother nor the convergence proof can be applied.} $\Box$
\end{remark}

A convergence proof of a variant of FAS for a class of monotone nonlinear elliptic problems is given by Hackbusch in~\cite{Hackbusch1985} and Reusken in~\cite{Reusken1988}. They proved convergence by linearising the FAS iteration and used the convergence theory for linear two-grid methods for symmetric elliptic problems as in~\cite{Bank1985}. Their proof was rigorous but requiring restrictive assumptions (the initial guess is close enough to the solution). Tai and Xu in~\cite{Tai2003,TaiXu2001} gave some uniform convergence estimates for a class of subspace correction methods applied to some nonlinear unconstrained and constraint convex optimization problems. But their methods is built upon nested finite element spaces and slightly expensive than FAS. Yavneh and Dardyk in~\cite{Yavneh2006} employed a simplified scalar analogy to provide an insight to the reason why FAS works but a rigorous proof is lacking. None of these theoretical work can be applied directly to our problem. We are investigating the convergence theory of FAS in different perspectives and will report our finding somewhere else.

\section{Numerical Experiments}
\label{SEC::NES}
In this section, some numerical results are presented to illustrate the efficiency of our multigrid method for the Darcy-Forchheimer model \eqref{Eqn::DF equation of the Prob}-\eqref{Eqn::N boundary of the Pro}. The following test problems are taken from~\cite{Jose2009}.
All of our experiments are implemented based on the $\rm{MATLAB}^{\copyright}$ software package $i$FEM~\cite{iFEM}. They were run on a laptop with a Inter i7-4720HQ 2.60GHz CPU and 16.0GB RAM.

We choose $\mu = 1$, $ \rho = 1$, $ \boldsymbol{K} = \boldsymbol{I}$, and $\Omega \subset \mathbb{R}^{\mathrm{2}}$ as the square $(-1,1)^2$. \LC{We use the uniform triangulation of $\Omega$.}

\begin{itemize}
\item Problem 1:
\begin{eqnarray*}
\boldsymbol{u}\left({x,y}\right) &=& \left[{x+y,x-y}\right]^{T},\\
p\left({x,y}\right) &=& x^3 + y^3,
\end{eqnarray*}

\begin{equation*}
\boldsymbol{f}\left( {x,y} \right) = \begin{bmatrix}
{\left( {1 + \beta \sqrt {2{x^2} + 2{y^2}} } \right)\left( {x + y} \right) + 3{x^2}}\\
{\left( {1 + \beta \sqrt {2{x^2} + 2{y^2}} } \right)\left( {x - y} \right) + 3{y^2}}
\end{bmatrix},
\end{equation*}
\begin{equation*}
g_N\left({x,y}\right) = \begin{cases}
 1 + y,  & \quad x = 1,\\ 
 1 - y,  & \quad x = -1,\\
 x - 1,  & \quad y = 1, \\
 -x - 1, & \quad y = -1.\\
 \end{cases} 
\end{equation*}

\item Problem 2:
\begin{eqnarray*}
\boldsymbol{u}\left({x,y}\right) &=& \left[{\frac{\left(x + 1 \right)^2}{4},-\frac{\left(x + 1 \right)\left(y + 1\right)}{2}}\right]^{T},\\
p\left({x,y}\right) &=& x^3 + y^3,
\end{eqnarray*}

\begin{equation*}
\boldsymbol{f}\left( {x,y} \right) = \begin{bmatrix}
{\quad \quad \,\,\frac{{{{\left( {x + 1} \right)}^2}}}{4}\left( {1 + \beta \frac{{\left( {x + 1} \right)}}{4}\sqrt {{{\left( {x + 1} \right)}^2} + 4{{\left( {y + 1} \right)}^2}} } \right) + 3{x^2}}\\
{ - \frac{{\left( {x + 1} \right)\left( {y + 1} \right)}}{2}\left( {1 + \beta \frac{{\left( {x + 1} \right)}}{4}\sqrt {{{\left( {x + 1} \right)}^2} + 4{{\left( {y + 1} \right)}^2}} } \right) + 3{y^2}}
\end{bmatrix},
\end{equation*}
\begin{equation*}
g_N\left({x,y}\right) = \begin{cases}
 1,  & \quad x = 1,\\ 
 0,  & \quad x = -1,\\
 -x - 1,  & \quad y = 1, \\
 0, & \quad y = -1.\\
 \end{cases} 
\end{equation*}
\end{itemize}
\LC{Numerically Problem 2 is harder to solve. Probably it is due to the fact that the initial guess, which is obtained by solving a linear Darcy system, is further away from the true solution.}

For all above test problems, $ g = 0 $. The chosen termination criterion is $$ r = r_{u} + r_{p} \leq tol ,$$ 
where 
\begin{eqnarray*}
r_{u} &= &\begin{cases}
 \left\| \boldsymbol{f} - \frac{\mu }{\rho }{\boldsymbol{K}^{ - 1}}\boldsymbol{u}_h^n + \frac{\beta }{\rho }\left| \boldsymbol{u}_h^n \right|\boldsymbol{u}_h^n + \nabla p_h^n \right\|/\left\| \boldsymbol{f} \right\|, &\quad\text{ when } \left\|\boldsymbol{f} \right\| \neq 0,\\
 \left\|\boldsymbol{f} - \frac{\mu }{\rho }{\boldsymbol{K}^{ - 1}}\boldsymbol{u}_h^n + \frac{\beta }{\rho }\left| \boldsymbol{u}_h^n \right|\boldsymbol{u}_h^n + \nabla p_h^n \right\|, &\quad\text{ when } \left\| \boldsymbol{f} \right\| = 0.
 \end{cases}\\
   r_{p} &=&\begin{cases} 
  \left\| g - \text{div}\boldsymbol{u}_h^n \right\|/\left\| g \right\|, &\quad\text{ when } \left\| g \right\|\neq 0,\\
  \left\| g - \text{div}\boldsymbol{u}_h^n \right\|, &\quad\text{ when } \left\| g \right\| = 0.\\
\end{cases} 
\end{eqnarray*}
 
\begin{figure}
  \centering
  \subfigure[Convergence rate by using multigrid solver]{
    \label{fig:subfig:data1beta30mgrate} 
    \includegraphics[width=2.3in]{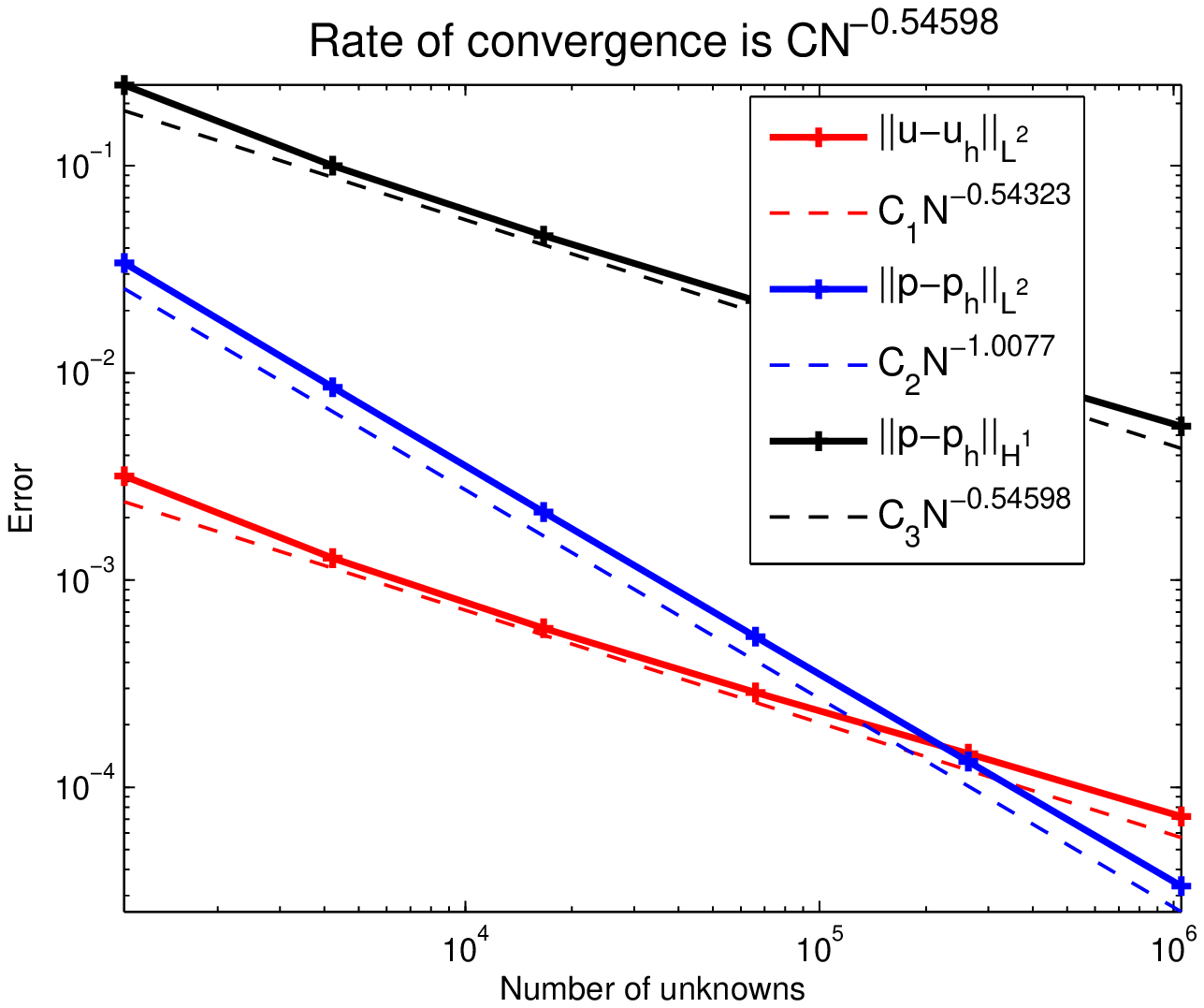}}
  \subfigure[Time complexity]{
    \label{fig:subfig:data1beta30timecomplexity} 
    \includegraphics[width=2.3in]{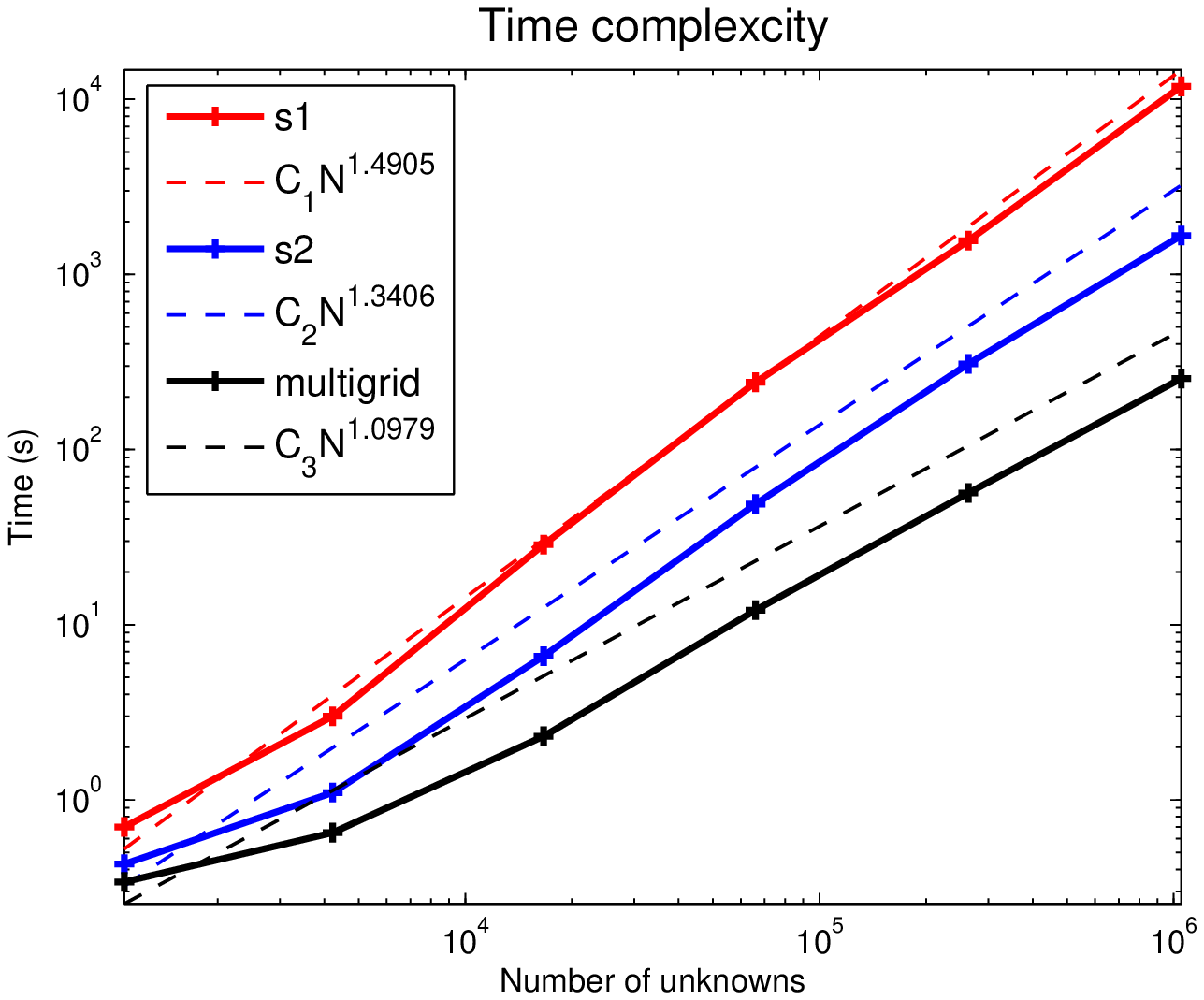}}
  \caption{Convergence rate by using multigrid solver and time complexity by using different solvers for Problem 1 with $\beta = 30$.}
  \label{fig:data1beta30rate} 
\end{figure}

\begin{figure}
  \centering
  \subfigure[Convergence rate by using multigrid solver]{
    \label{fig:subfig:data2beta30mgrate} 
    \includegraphics[width=2.3in]{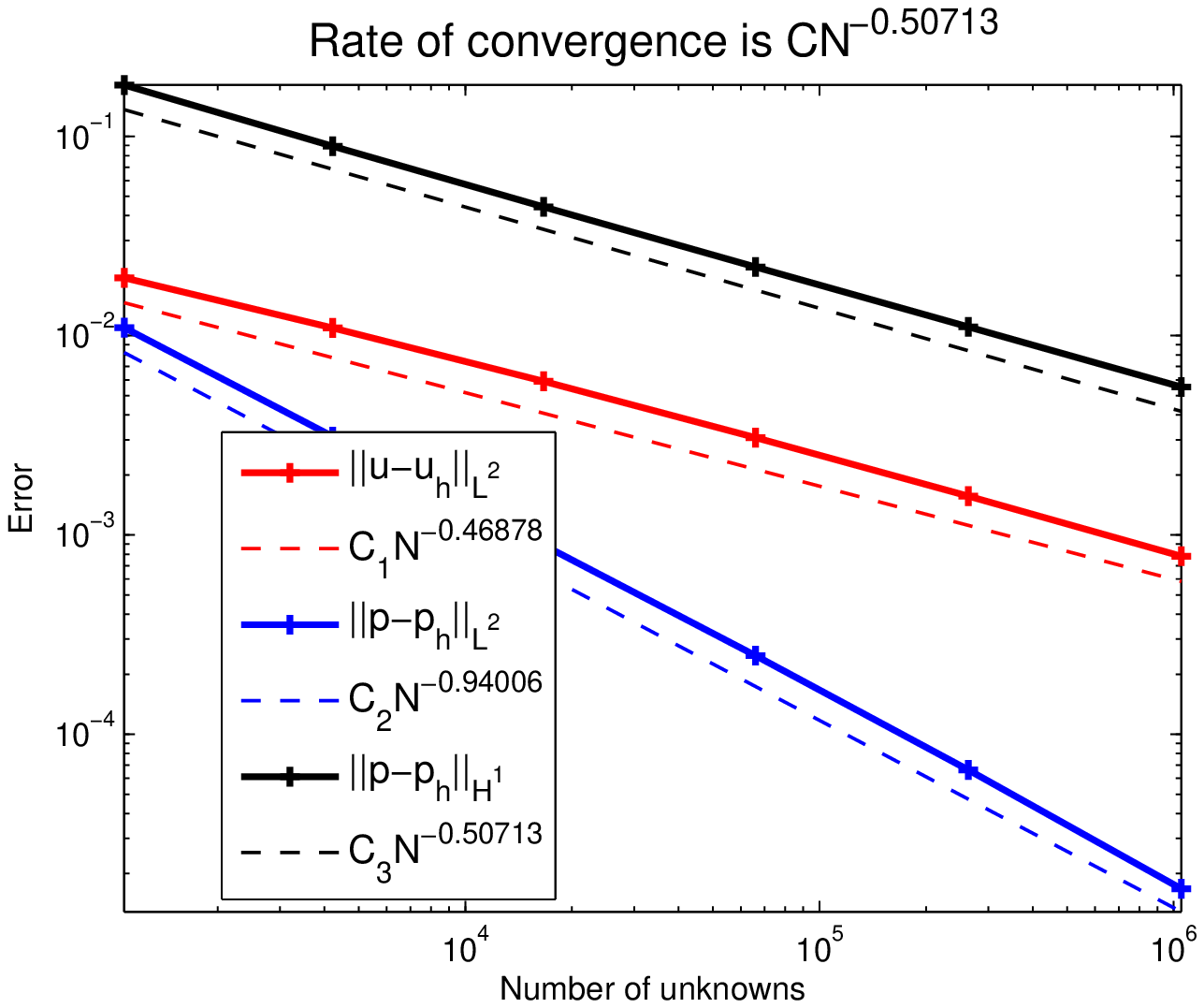}}
  \subfigure[Time complexity]{
    \label{fig:subfig:data2beta30timecomplexity} 
    \includegraphics[width=2.3in]{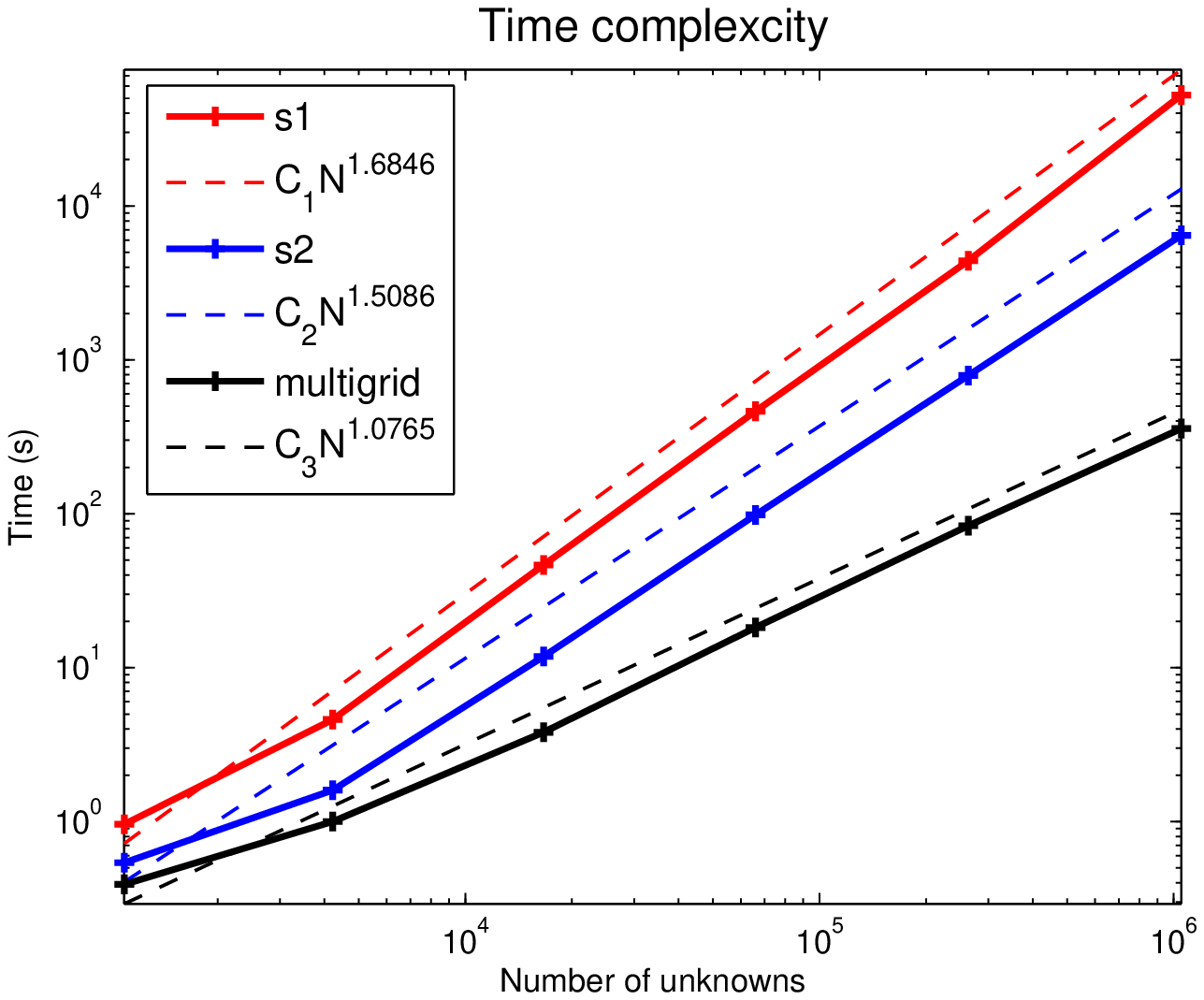}}
  \caption{Convergence rate by using multigrid solver and time complexity by using different solvers for Problem 2 with $\beta = 30$.}
  \label{fig:data2beta30rate} 
\end{figure}

\LC{We first use the accuracy test to confirm that our nonlinear multigrid iteration will convergent to an approximation of the problem of consideration.}
In the following experiments, the letter $ N $ stands for `Number of unknowns of $ p $', which is the same as `Numbers of vertices', so $ h = \frac{2}{\sqrt{N} - 1 } $, which represents the discretization mesh size in one direction. Numerical results, see Fig.~\ref{fig:subfig:data1beta30mgrate},~\ref{fig:subfig:data2beta30mgrate}, confirmed the convergence order for $ {\left\| {\boldsymbol{u} - \boldsymbol{u}_{h}} \right\|_{{L^2}}} $ and $ {\left\| {p - {p_h}} \right\|_{{H^1}}} $ are $O\left( h \right) = O(N^{1/2})$.  
The accuracy of the pressure approximations, however, is not as good as that of velocity. Meanwhile, in consideration of the computation cost, the sufficiently accurate results were achieved when $ tol = 10^{-6} $ for Problem 1 and 2. \LC{The stopping tolerance can be varying in different levels to further reduce the cost. A guide line is below the truncation error~\cite{Brandt1977}.} The authors in~\cite{Jose2009}, however, use $ tol = 1.95h $, which is only enough for the $ L^{2} $-norm approximation for velocity. We shall use $ tol = 10^{-6}$ in the remaining numerical experiments. 
 
\LC{ For all tests, the iteration steps and CPU time of each solver are listed in tables. We are aware that the CPU time depends on the implementation and testing environment: the programming language, optimization of codes, and the hardware (memory and cache), etc. Our code has been optimized using vectorization technique and all results were measured and compared in the same test environment so that the CPU time could be a good indicator of the efficiency. The CPU time will be also used to find the asymptotic time complexity of each method; see Fig.~\ref{fig:subfig:data1beta30timecomplexity},~\ref{fig:subfig:data2beta30timecomplexity}.}

As it is \LC{proved} in~\cite{Wheeler2008}, the PR nonlinear iteration converges for any $\alpha > 0 $. Its rate of convergence, however, is very sensitive to the choice of this parameter. From the convergence proof of the PR iteration in~\cite{Wheeler2008}, we inferred that the choices of $ \alpha $ depends on \JH{the Forchheimer number $ \beta $ which controls the magnitude of the nonlinearity as $\rho$ is fixed}. We give an empirical choice of parameter $ \alpha = 1/\beta $ and compared with the choice $\alpha = 1$ suggested in \cite{Jose2009} in Table \ref{TAB::Comparison alpha beta 10-30} and \ref{TAB::Comparison alpha beta 40-60}. As shown in Table \ref{TAB::Comparison alpha beta 10-30} and \ref{TAB::Comparison alpha beta 40-60}, this choice of the parameter $ \alpha $ is much better than the fixed selection for different values of $ \beta $. Therefore, this choice of $ \alpha $ will be used in the remaining numerical experiments. 

\begin{table}
\caption{Comparison of different values of $\alpha$ in PR iteration with $h = \frac{1}{64}$ for $ \beta = 10, 20, 30$.}
\label{TAB::Comparison alpha beta 10-30}
\begin{tabular}{*{10}{l}}
\hline\hline\noalign{\smallskip}
\multirow{2}*{Problem} & &\multicolumn{2}{ c }{$ \beta = 10 $}&\multicolumn{2}{ c }{$ \beta = 20 $}&\multicolumn{2}{ c }{$ \beta = 30 $}\\
                        &                    & $ \alpha = 1 $ & $ \alpha = 1/10 $ & $ \alpha = 1 $ & $ \alpha = 1/20 $ & $ \alpha = 1 $ & $ \alpha = 1/30 $ \\
                        \noalign{\smallskip}\hline\noalign{\smallskip}
 
\multirow{2}*{Problem 1}& iter         & 229       & 73           & 457       & 105    & 686       & 120 \\
                        & CPU time     & 14 s      & 4 s          & 26 s      & 6 s    & 38 s      & 7 s \\
                        \hline\noalign{\smallskip}
\multirow{2}*{Problem 2}&iter          & 230       & 171          & 459       & 183    & 688       & 191 \\
                        &CPU time      & 13 s      & 10 s         & 26 s      & 11 s   & 38 s      & 11 s \\   
\noalign{\smallskip}\hline\hline
\end{tabular}
\end{table}

\begin{table}
\caption{Comparison of different values of $\alpha$ in PR iteration with $h = \frac{1}{64}$ for $ \beta = 40, 50, 60$.}
\label{TAB::Comparison alpha beta 40-60}

\begin{tabular}{*{10}{l}}
\hline\hline\noalign{\smallskip}

\multirow{2}*{Problem} & &\multicolumn{2}{ c }{$ \beta = 40 $}&\multicolumn{2}{ c }{$ \beta = 50 $}&\multicolumn{2}{ c }{$ \beta = 60 $}\\ 
                        &                    & $ \alpha = 1 $ & $ \alpha = 1/40 $ & $ \alpha = 1 $ & $ \alpha = 1/50 $ & $ \alpha = 1 $ & $ \alpha = 1/60 $ \\
                        \noalign{\smallskip}\hline\noalign{\smallskip} 
\multirow{2}*{Problem 1}& iter         & 914       & 126          & 1143      & 129       & 1371       & 131 \\
                        & CPU time     & 53 s      & 7 s          & 66 s      & 7 s       & 79 s       & 8 s \\
                        \hline\noalign{\smallskip}
\multirow{2}*{Problem 2}&iter          & 917       & 198          & 1146      & 205       & 1376       & 213 \\
                        &CPU time      & 52 s      & 11 s         & 65 s      & 11 s      & 79 s       & 12 s \\
                        \noalign{\smallskip}\hline\hline  
\end{tabular}
\end{table}

We then compare the FAS multigrid method using PR as smoother with the PR iterative method for solving this nonlinear system. Here we choose $ m = 3 $ for all the following tests. It means that we apply three PR iterations in the pre-smoothing step and post-smoothing step, respectively. \LC{Each V-cycle step is approximately $9$ PR iterations ($6$ for the finest level and $3$ for iterations in all coarser levels as the size of the system is reduced by $1/4$) in terms of complexity.} In order to keep the symmetry of the V-cycle, we switch the ordering of the linear and nonlinear steps of the PR iteration in the post-smoothing step. We set $ h =1/16 $ as the coarsest mesh and solve the nonlinear problem in the coarsest mesh using PR iteration. 

The PR solver is denoted by \JH{pr}, whereas the multigrid solver is denoted by \LC{mg}. I - number of iterations, and CPU - CPU time. `s1' represents that we solve these linear saddle point systems \eqref{Eqn::matrix form of nonlinear iteration} directly in each step, `s2' is that we solve the primal SPD system \eqref{Eqn::SPD system} mentioned in Section \ref{SEC::Nonlinear iteration} rather than solving the saddle point system. `mg' stands for our multigrid solver, in which the PR iteration is constructed based on `s2'. In all examples we achieve optimal order convergence of $ {\left\| {\boldsymbol{u} - \boldsymbol{u}_{h}} \right\|_{{L^2}}} $ and $ {\left\| {p - {p_h}} \right\|_{{H^1}}}$. Compared with the PR iteration, we can obtain the same accuracy by using our multigrid method with less iterations. We can get similar results for different values of the Forchheimer number $\beta$.

Since our focus is on the efficiency of solvers, we mainly report the comparison of the number of iterations and CPU time by using different solvers. Numerical tests were performed for several cases of different values of the Forchheimer number $\beta$ for Problem 1 and 2, and the behavior of these experiments is similar for all chosen cases. \JH{All problems are becoming harder to solve as the Forchheimer number $\beta$ increases, mainly because $\beta $ enhances the nonlinearity.} Therefore, without loss of critical substance and clarity, here we only show the results for $ \beta = 30 $ to demonstrate the merits of our method.

It can be observed that our multigrid solver required significantly fewer iterations and CPU time than the other two solvers in Table \ref{TAB::CPU time of P1 beta30} and \ref{TAB::CPU time of P2 beta30}. \JH{More importantly, iteration steps are uniformly stable with respect to $h$} and the time complexity of our multigrid solver is nearly linear, i.e., $ O(N) $, shown in Fig. \ref{fig:data1beta30rate} and \ref{fig:data2beta30rate}. In contrast, for the PR methods, \JH{iteration steps increase as $h$ decreases} and the time complexity seems to be more than linear. \JH{For the largest size we have tested, our multigrid solver is more than $40$ times faster than the original PR iteration.} In Table \ref{TAB::mg iteration steps of P1} and \ref{TAB::mg iteration steps of P2}, the number of iterations are compared for different values of $\beta$ and it is demonstrated that our multigrid method is also robust to both mesh size $h$ and the Forchheimer number $\beta$ while \LC{PR iteration is not, see Table 1 and 2.} \JH{It is worth noting that even for a linear Stokes type equation, construct a solver robust to a critical parameter is not easy~\cite{Mardal2004,Olshanskii2006}.} 

\begin{table}
\caption{Comparison of number of iterations and CPU time of Problem 1 by using different solvers with $\beta = 30$.}
\label{TAB::CPU time of P1 beta30}
\begin{tabular}{lllllll} 
\hline\hline\noalign{\smallskip}
$ h $             & DoFs     & I(pr) & I(mg) & CPU(s1)    & CPU(s2)     & CPU(mg) \\ 
\noalign{\smallskip}\hline\noalign{\smallskip}

$\frac{1}{16} $   & 5,185    &  50   &  1    &  0.70 s    &  0.43 s     & 0.34 s \\

$\frac{1}{32} $   & 20,609   &  81   &  6    &  3.0 s     &  1.1 s      & 0.65 s \\

$\frac{1}{64} $   & 83,177   &  120  &  6    &  28.6 s    &  6.6 s      & 2.3 s \\

$\frac{1}{128}$ & 328,193    &  154  &  6    &  242.3 s   &  48.8 s     & 12.1 s \\

$\frac{1}{256}$ & 1,311,745  &  168  &  6    &  1554.7 s  &  308.3 s    & 56.5 s \\
 
$\frac{1}{512} $ & 5,244,929 &  185  &  5    &  11857.3 s &  1667.7 s   &  254.6 s \\
\noalign{\smallskip}\hline\hline
\end{tabular}
\end{table}

\begin{table}
\caption{Comparison of iteration steps of multigrid solver according to different $h$ and $\beta$ for Problem 1 with $ \alpha = 1/\beta $.}
\label{TAB::mg iteration steps of P1}
\begin{tabular}{llllll} 
\hline\hline\noalign{\smallskip}
$ h $             & $\beta = 10$& $\beta = 20$ & $\beta = 30$ & $\beta = 40$ & $\beta = 50$ \\
 \noalign{\smallskip}\hline\noalign{\smallskip}

$\frac{1}{32} $   & 4           &  6           &  6         &  7             &  7  \\

$\frac{1}{64} $   & 4           &  6           & 6          &  7             &  7  \\

$\frac{1}{128}$   & 4           &  5           & 6          &  6             & 7   \\

$\frac{1}{256}$   & 4           &  5           & 6          & 6              & 6  \\

$\frac{1}{512} $  & 3           &  5           &  5         & 6              & 6  \\
\noalign{\smallskip}\hline\hline
\end{tabular}
\end{table}

\begin{table}
\caption{Comparison of number of iterations and CPU time of Problem 2 by using different solvers with $\beta = 30$.}
\label{TAB::CPU time of P2 beta30}
\begin{tabular}{lllllll} 
\hline\hline\noalign{\smallskip}
$ h $             & DoFs     & I(pr) & I(mg)& CPU(s1)       & CPU(s2)     & CPU(mg) \\ 
\noalign{\smallskip}\hline\noalign{\smallskip}

$\frac{1}{16} $   & 5,185    &  92   &  1   &  0.96 s       &  0.54 s     & 0.39 s  \\
 
$\frac{1}{32} $   & 20,609   &  128  &  9   &  4.6 s        &  1.6 s      & 1.0 s \\

$\frac{1}{64} $   & 83,177   &  191  &  9   &  46.5 s       &  11.8 s     & 3.8 s\\

$\frac{1}{128}$ & 328,193    &  296  &  9   &  462.9 s      &  98.3 s     & 18.2 s\\

$\frac{1}{256}$ & 1,311,745  &  468  &  8   &  4412.9 s     &  792.6 s    & 83.6 s\\

$\frac{1}{512} $ & 5,244,929 &  746  &  7   & $> 14$ hours  &  6440.3 s   &  357.2 s \\
\noalign{\smallskip}\hline\hline
\end{tabular}
\end{table}

\begin{table}
\caption{Comparison of iteration steps of multigrid solver according to different $h$ and $\beta$ for Problem 2 with $ \alpha = 1/\beta $.}
\label{TAB::mg iteration steps of P2}
\begin{tabular}{llllll} 
\hline\hline\noalign{\smallskip}
$ h $             & $\beta = 10$& $\beta = 20$ & $\beta = 30$ & $\beta = 40$ & $\beta = 50$ \\ 
\noalign{\smallskip}\hline\noalign{\smallskip}
$\frac{1}{32} $   & 5           &  7           &  9         &  11            &  12  \\

$\frac{1}{64} $   & 5           &  7           & 9          &  11            &  12  \\

$\frac{1}{128}$   & 5           &  7           & 9          &  10            & 11   \\

$\frac{1}{256}$   & 4           &  6           & 8          &  9             & 10  \\

$\frac{1}{512} $  & 4           &  5           &  7         &  8             & 9  \\
\noalign{\smallskip}\hline\hline
\end{tabular}
\end{table}

\section{Conclusions}
\label{SEC:Conclusions}
In this paper, we constructed a nonlinear multigrid method for a mixed finite element method of the two-dimensional Darcy-Forchheimer model. We presented a comparative study between the multigrid solver and the PR iterative solver, at the same time compared CPU time of the efficient solver of solving the SPD systems with that obtained by solving the linear saddle point systems directly. We took into account the pressure accuracy when we set the termination criterion, and chose a better value of the stopping criterion $tol$. In comparison with the authors in~\cite{Jose2009} always chose $ \alpha = 1$ for different values of the Forchheimer number $ \beta $, we reported a better choice and compared with the previous choice through comparing the number of iterations and CPU time. The results obtained from our tests indicate that the multigrid solver is very efficient for numerically solving this nonlinear elliptic equation. The number of iterations and CPU time for using multigrid solver are shown to be significantly less than that obtained by using \JH{the PR iteration alone}.

In the future work, we shall extend our results to \JH{three directions}. One is that we would like to find a better smoother, which is used in the pre-smoothing and post-smoothing step, to further reduce CPU time and make the multigrid solver more efficient. Another is that we intend to carry out some studies on the three-dimensional Darcy-Forchheimer problem and the real application in a porous medium. \JH{We shall also investigate the theoretical study of the convergence proof of FAS.}




\end{document}